\documentclass[12pt]{article}
\usepackage{group21, times}
\usepackage{amstex}
\usepackage{theorem}
\usepackage{amssymb}
\usepackage{amscd}

\title{The idea of locality}
\author{Victor G. Kac\thanks
{Partially supported by NSF grant DMS-9622870}}
\address{Department of Mathematics, MIT,\\ 
Cambridge, MA 02139}

\newcommand{\ccc}{{\Bbb C}}
\newcommand{\zz}{{\Bbb Z}}
\newcommand{\jj}{{\mathfrak{g}}}
\newcommand{\jjb}{\tilde{{\mathfrak{g}}}}

\newcommand{\pa}{\partial}
\newcommand{\de}{\delta}
\newcommand{\dd}{{\cal{D}}}
\newcommand{\ff}{F}
\newcommand{\rr}{R}

\newcommand{\Res}{\operatorname{Res}}
\newcommand{\Vir}{\operatorname{Vir}}
\newcommand{\ddiv}{\operatorname{div}}
\newcommand{\ddim}{\operatorname{dim}}
\newcommand{\ddeg}{\operatorname{deg}}

\newcommand{\ind}{}

\newcommand{\pr}[1]{\ind{\it #1.\/}}
\newcommand{\qed}{\hfill$\Box$}
\newcommand{\Qed}{\qed\smallskip}

\newtheorem{proposition}{Proposition}[section]
\newtheorem{lemma}{Lemma}[section]

\newtheorem{theorem}{Theorem}[section]
\newtheorem{cor}{Corollary}[section]

{\theorembodyfont{\rmfamily} 
  \newtheorem{example}{Example}[section]
  \newtheorem{remark}{Remark}[section]
}

\newenvironment{examples}{\bigskip\noindent {\bfseries
    Examples.}\enspace}{\medskip}

\renewcommand{\geq}{\geqslant}
\renewcommand{\leq}{\leqslant}
\newcommand{\gcn}{gc_N}
\newcommand{\rrt}{\tilde{R}}
\newcommand{\nt}{\tilde{n}}

\newcommand{\alphalist}{% changes enumerate 1st level to a)...z)
  \renewcommand{\theenumi}{\rm\alph{enumi}}%
  \renewcommand{\labelenumi}{\theenumi)}%
}

\newcommand{\alphaparenlist}{% changes enumerate 1st level to (a)...(z)
  \renewcommand{\theenumi}{\rm\alph{enumi}}%
  \renewcommand{\labelenumi}{\theenumi)}%
}

\newcounter{bean}
\newenvironment{deflist}[1]%
    {
    \begin{list}{\bf #1\arabic{bean}}
        {\usecounter{bean}
          \labelsep=1em
          \settowidth{\labelwidth}{#1\thebean:}
          \addtolength{\labelwidth}{1.1ex} 
          \leftmargin=\labelwidth 
          \addtolength{\leftmargin}{\labelsep} }

}
    {\end{list}}

\newcommand{\req}[1]{(\ref{#1})}

\begin{document}
\maketitle
\begin{abstract}
This is a review of recent results on conformal (super)algebras. It may be
viewed as an amplification of my Wigner medal acceptance speech (given in
July 1996 in Goslar, Germany) reproduced in the introduction.
\end{abstract}

\setcounter{section}{-1}
\section{Introduction.}
\label{sec:0}

It is a well kept secret that the theory of Kac-Moody algebras has
been a disaster.

True, it is a generalization of a very important object,
the simple finite-dimen\-sional Lie algebras, but a generalization too
straightforward to expect anything interesting from it.  True, it is
remarkable how far one can go with all these $e_{i}$'s, $f_{i}$'s and
$h_{i}$'s.  Practically all, even most difficult results of
finite-dimensional theory, such as the theory of characters, Schubert
calculus and cohomology theory, have been extended to the general
set-up of Kac-Moody algebras.  But the answer to the most important
question is missing:  what are these algebras good for?  Even the most
sophisticated results, like the connections to the theory of quivers,
seem to be just scratching the surface. 

However, there are two notable
exceptions.  The best known one is, of course, the theory of affine
Kac-Moody algebras.  This part of the Kac-Moody theory has
deeply penetrated many branches of mathematics and physics.  The most
important single reason for this success is undoubtedly the
isomorphism of affine algebras and central extensions of loop
algebras, often called current algebras.

The second notable exception is provided by Borcherds' algebras which
are roughly speaking the spaces of physical states of certain chiral
algebras.

At this point a natural question arises:  what do these notable
exceptions have in common?  The answer, in my opinion, lies in the
idea of locality.

The concept of locality is a beautiful synthesis of three ideas:

\begin{enumerate}
\item Einstein's special relativity postulate,
\item Heisenberg's uncertainty principle,
\item the notion of quantum field.
\end{enumerate}

As far as I know, the concept of locality was first rigorously
formulated by Pauli and then incorporated by Wightman in his
axiomatics of quantum field theory.

The locality axiom states that quantum fields, whose supports are
space-like separated, must commute.

The locality axiom is empty for $1$-dimensional space-time, but
becomes a very complicated condition for higher-dimensional space-times.  
However, in the $2$-dimen\-sional case the light cone is a union
of two straight lines and one may consider chiral quantum fields, that
is, fields depending on one null coordinate.

For chiral fields the locality axiom reduces to a very simple
algebraic condition.  It just states that the commutator of chiral
fields is a finite linear combination, with coefficients being chiral
fields, of the delta-function and its derivatives.

Coming back to our two notable exceptions, we see that both are
spanned by Fourier coefficients of pairwise local chiral fields!

Now, if we look around, we will immediately find that all the
important Lie algebras (or superalgebras) satisfy this locality
property:  the Virasoro algebra and its marvelous super
generalizations, the $W$-algebras, the Lie algebra $W_{1+\infty}$,
etc., etc.

This leads us to the program of study of Lie algebras and
superalgebras satisfying the sole locality property.  It is an
extensive program, but I hope that its main features will be worked
out by the year 2000.

The basic message I wanted to convey by my speech is this.  Some of
the best ideas come to my field from the physicists.  And on top
of this they award me a medal.  One couldn't hope for a better deal.

Thank you.

\section{Preliminaries on Lie superalgebras of formal distributions
and conformal superalgebras.}
\label{sec:1}

A {\it formal distribution} 
(usually called a field by physicists) with coefficients 
in a complex vector space $U$ is a generating series
of the form:
\begin{equation*}
a(z)=\sum_{n\in\zz}a_{(n)}z^{-n-1},
\end{equation*}
where $a_{(n)}\in U$ and $z$ is an
indeterminate.

Two formal distributions $a(z)$ and $b(z)$ with coefficients in a Lie
superalgebra $\jj$ are called (mutually) 
{\it local} if for some $N\in \zz_{+}$ one has: 
\begin{equation}
(z-w)^N[a(z),b(w)]=0 \text{ for some } N\in\zz_+ \, .
\label{1.1}
\end{equation}
Introducing the {\it formal delta-function}
\begin{equation*}
\de(z-w)=z^{-1}\sum_{n\in\zz}\left(\frac{z}{w}\right)^n \, ,
\end{equation*}
we may write a condition equivalent to (\ref{1.1}):
\begin{equation}
[a(z),b(w)]=\sum_{j=0}^{N-1}\left(a_{(j)}b\right)(w)\pa^j_w\de(z-w)/j!
\label{1.2}
\end{equation}
for some formal distributions $(a_{(j)}b)(w)$ (\cite{1}, Theorem 2.3),
which are uniquely determined by the formula
\begin{equation}
\left(a_{(j)}b\right)(w)=\Res_{z}(z-w)^j[a(z),b(w)] \, .
\label{1.3}
\end{equation}
Formula (\ref{1.3}) defines a $\ccc$-bilinear product $a_{(j)}b$ for each 
$j\in\zz_{+}$ on the space of all formal distributions with coefficients in 
$\jj$.
 
Note also that the space (over $\ccc$) of all formal distributions with 
coefficients in $\jj$ is a (left) module over $\ccc[\pa]$, where the
action of $\pa = \pa_z$ is defined in the obvious way,
so that $\pa_z a(z)=\sum_n(\pa a)_{(n)}z^{-n-1}$, where
$(\pa a)_{(n)}=-na_{(n-1)}$.

The Lie superalgebra $\jj$ is called a {\it Lie superalgebra of formal 
distributions} if there exists a family $\ff$ of pairwise local formal 
distributions whose coefficients span $\jj$. In such a case we say that
the family $\ff$ {\it spans} $\jj$.
We will write $(\jj,\ff)$ to emphasize the dependence of
$\ff$.  

The simplest example of a Lie superalgebra of formal distributions is the
{\it current superalgebra} $\jjb$ associated to a Lie superalgebra 
$\jj$:
\begin{equation*}
\jjb=\ccc\left[t,t^{-1}\right]\otimes\jj \, .
\end{equation*}
It is spanned by the following family of pairwise local formal distributions
$(a\in {\jj})$:  
\begin{equation*}
a(z)=\sum_{n\in\zz}\left(t^n\otimes a\right)z^{-n-1}.
\end{equation*}
Indeed, it is immediate to check that
\begin{equation*}
[a(z),b(w)]=[a,b](w)\de(z-w).
\end{equation*}

The simplest example beyond current algebras is the (centerless) {\it Virasoro
algebra}, the Lie algebra with the basis $L_{n}$ $(n\in\zz)$ and 
commutation relations
\begin{equation*}
\left[ L_m,L_n\right]=(m-n)L_{m+n} \, .
\end{equation*}
It is spanned by the local formal distribution 
$L(z)=\sum_{n\in\zz}L_n z^{-n-2}$, since one has:
\begin{equation}
[L(z),L(w)]=\pa_w L(w)\de(z-w)+2L(w)\pa_w\de(z-w) \, .
\label{1.4}
\end{equation}

Note that $\ccc[\pa]\ff$ is a $\ccc[\pa]$-submodule of the space of formal 
distributions which still consists of pairwise local formal distributions.
The Lie superalgebra $\jj$ of formal distributions spanned by $\ff$ is called
{\it simple} if $\jj$ contains no non-trivial ideals spanned by a subspace
of $\ccc[\pa]\ff$.

The Virasoro algebra has no ideals at all, but the current algebra always has
``evaluation'' ideals. Nevertheless the current algebra $\jjb$ associated to 
Lie superalgebra $\jj$ is simple in the above sense iff $\jj$ is
simple in the usual sense.

Given a Lie superalgebra of formal distributions $(\jj,\ff)$, we may always 
include $\ff$ in the minimal family $\ff^c$ of pairwise local formal 
distributions which is closed under $\pa$ and all products
(\ref{1.3}) (\cite{1}, \S~2.7).

We say that a Lie superalgebra $(\jj,\ff)$ of 
formal 
distribution is  $finite$ if $\ff^c$  is a finitely generated 
$\ccc[\pa]$-module.

The Virasoro algebra and the 
current (super)algebras associated to finite-dimen\-sional Lie (super)algebras 
are 
finite. The finiteness
condition provides the choice of the ``non-twisted'' moding, as the following
example shows.

%%%\examp{Example 1.1.} \  
\begin{example}
\label{ex:1.1}
Let $\jj$ be a finite-dimensional Lie algebra
and  let $\jj=\jj^0+\jj^1$ be a 
$\zz/2\zz$-gradation.
Then
\begin{equation*}
\jjb:=\ccc\left[t,t^{-1}\right]\otimes\jj^0+t^{\frac{1}{2}}
\ccc\left[t,t^{-1}\right]\otimes
\jj^1
\end{equation*}
is a subalgebra of the Lie algebra $\ccc\left[t^{\frac{1}{2}},
t^{-\frac{1}{2}}\right]\otimes\jj^0$. This is a ``twisted'' current algebra. 
It is spanned by the set $\ff$ of pairwise local formal distributions 
\begin{equation*}
%\begin{split}
a(z)=
%&
\sum_{n\in\zz}\left(t^n\otimes a\right)z^{-n-1} \text{ for } 
a\in\jj^0, \quad \text{and }
%\\
a(z)=
%&
\sum_{n\in\zz}\left(t^{n+\frac{1}{2}}\otimes a\right)z^{-n-1} 
\text{ for } a\in \jj^1 \, .
%\end{split}
\end{equation*} 
However, $[a(z),b(w)]=w[a,b](w)\de(z-w)$ if $a,b\in\jj^1$. 
Hence $\jj$ is not a finite Lie algebra of formal distributions, at least
if $\jj=[\jj,\jj]$, since $\ff^c$ contains all formal distributions 
$w^n a(w)$, $n\in\zz_+$.
\end{example}

Lie superalgebras of formal distributions can be studied via 
conformal superalgebras introduced in \cite{1}. 

A {\it conformal superalgebra} $\rr$ is a left 
($\zz/2\zz$-graded) $\ccc[\pa]$-module 
$\rr$ with
a  $\ccc$-bilinear product $a_{(n)}b$ for each $n\in\zz_{+}$ such that the 
following axioms hold ($a,b,c\in\rr$, $m,n\in\zz_{+}$):
%%%\begin{equation*}
%%%\begin{split}
\begin{deflist}{C}
\setcounter{bean}{-1}
\item %%%\text{(C0)}&\\
  $a_{(n)}b=0 \text{ for } n \gg 0$,
\item
  %%% \text{(C1)}&
  $(\pa a)_{(n)}b=-na_{(n-1)}b$,
\item 
  %%%\text{(C2)}&
  $a_{(n)}b=(-1)^{p(a)p(b)} \sum_{j=0}^{\infty}
  (-1)^{j+n+1}(\pa^j/j!)b_{(n+j)}a$,
\item
  %%%\text{(C3)}&
  $a_{(m)}\left(b_{(n)}c\right) = \sum_{j=0}^{\infty}\binom{m}{j}
  \left(a_{(j)}b\right)_{(m+n-j)}c+(-1)^{p(a)p(b)}b_{(n)}
  \left(a_{(m)}c\right)$.
\end{deflist}
%%%\end{split}
%%%\end{equation*}

Of course, a conformal algebra coincides with its even part, i.e. $p(a)=0$ 
for all $a\in\rr$ in this case. Note the following consequence of (C1) and
(C2):
%%%\begin{equation*}
%%%\text{(C1$'$)}
\begin{description}
\item[(C1')\hskip1ex] 
$a_{(n)}\pa b=\pa\left(a_{(n)}b\right)
+n a_{(n-1)}b$.
%%%\hspace{4.7cm}
\end{description}
%%%\end{equation*}
It follows that $\pa$ is a derivation of all products~(\ref{1.3}).

It is shown in \cite{1},Sec.2.7 that if $(\jj,\ff)$ is a Lie superalgebra
of formal distributions, then $\ff^c$ is a conformal superalgebra with
respect to these products.

Conversely, if $R=\oplus_{i\in I}\ccc[\pa]a^i$ is a free as 
$\ccc[\pa]$-module conformal superalgebra, 
we may associate to $R$ a Lie superalgebra
of formal distributions $(\jj(\rr),\ff)$ with the basis $a_{(m)}^i$
($i\in I$, $m\in\zz$)
and $\ff=\left\{a^i(z)=\sum_{n}a_{(n)}^i z^{-n-1}\right\}_{i\in I}$  
with the bracket (cf.~(\ref{1.2})):
\begin{equation}
[a^i(z),a^j(w)]=\sum_{k\in\zz_{+}}\left(a^i_{(k)} a^j\right)(w)
\pa_{w}^k\de(z-w)/k! \, ,
\label{1.5}
\end{equation}
so that $\ff^c=\rr$. 

It is well-known that a finitely generated $\ccc[\pa]$-module $\rr$ is a
direct sum of $r$ copies of $\ccc[\pa]$ and a $d$-dimensional (over $\ccc$)
$\ccc[\pa]$-invariant subspace. The number $r$ is called the {\it rank}
of $\rr$, and the pair of numbers $(r,d)$ is called the {\it size} of
$\rr$. A a factor module of $\rr$ by a 
non-zero submodule has smaller size $(r',d')$ in the sense that either
$r'<r$ or $r'=r$ but $d'<d$. We shall say that $\rr$ is {\it finite} if both
numbers $r$ and $d$ are finite (which is equivalent to $\rr$ being a
finitely generated $\ccc[\pa]$-module).

The {\it center} of $\rr$ is the set of the elements $a\in\rr$ such that
$a_{(n)}\rr=0$ for all $n\in\zz_+$ (then, by (C2), $\rr_{(n)}a=0$ for all
$n\in\zz_+$ as well).

It is shown in \cite{1}, Proposition 2.7 that a finite (as a 
$\ccc[\pa]$-module) conformal superalgebra $\rr$ with a trivial center is a 
free $\ccc[\pa]$-module (cf.\ Lemma 2.1 in \S\ref{sec:2}).
Hence the above discussion implies the following proposition.      

\begin{proposition}
\label{prop:1.1}
Every finite superalgebra of formal distributions $(\jj, \ff)$ with trivial
center is a quotient of a Lie superalgebra $\jj(\rr)$, where $\rr$ is a finite
conformal superalgebra with trivial center, by an ideal that does not contain
all $a_{(n)}$, $n\in\zz$, for a non-zero element $a\in\rr$. For such Lie
superalgebras the $\ccc[\pa]$-module $\ff^c$ is free.
\end{proposition}

The simplest example of a finite conformal superalgebra is the {\it current
conformal superalgebra} associated to a finite-dimensional Lie superalgebra
$\jj$:
\begin{equation*}
\rr(\jjb)=\ccc[\pa]\otimes_{\ccc}\jj \, ,
\end{equation*}
which has size $(\ddim\jj,0)$, with the products defined by:
\begin{equation*}
a_{(0)}b=[a,b] \, \quad
a_{(j)}b=0 \text{ for } j>0 \, , \quad 
a,b\in\jj \, ,
\end{equation*}
and the {\it Virasoro conformal algebra} $\Vir=\ccc[\pa]L$, which has size
$(1,0)$, with the products (cf.~(\ref{1.4})):
\begin{equation*}
L_{(0)}L=\pa L \, , \quad
L_{(1)}L=2L \, , \quad
L_{(j)}L=0 \text{ for } j>1 \, .
\end{equation*}
More complicated finite conformal superalgebras are constructed
in \S\ref{sec:5}.

Probably the most important infinite conformal algebra is that associated
to the Lie algebra $\dd$ of differential operators on the circle. Note also
that any vertex (=chiral) algebra \cite{2} may be viewed as a conformal 
superalgebra if we forget about the products $a_{(n)}b$ for $n<0$ and let
$\pa=T$, the infinitesimal translation operator (see \cite{1}, \cite{2}),
but this conformal superalgebra is always infinite.

Proposition 1.1 reduces the study of finite superalgebras of formal
distributions to that of finite conformal superalgebras. All notions 
concerning Lie superalgebras are automatically translated into the 
language of conformal superalgebras and vice versa.

A {\it subalgebra} $S$ (resp.\ {\it ideal} $I$) of a conformal 
superalgebra $\rr$
is a $\ccc[\pa]$-submodule of $R$ such that for $a,b\in S$ (resp.\ $a\in\rr$,
$b\in I$) we have $a_{(n)}b\in S$ (resp.~$\in I$) for all $n\in\zz_+$.
Note that 
the left (or right) ideal is automatically $2$-sided due to (C2).
The notions of a homomorphism and an isomorphism of conformal
superalgebras are obvious.A conformal superalgebra $\rr$ is called {\it simple}
if its only ideals are $0$ and $\rr$ and not all of the products are trivial.

A Lie superalgebra of formal distributions $(\jj,\ff)$ with a closed under
all products (\ref{1.3}) family $\ff$ is simple iff the associated conformal 
superalgebra $\ff^c$ is simple. For example, Vir is a simple conformal
algebra and the current algebra $\rr(\jjb)$ is a simple conformal 
{\it super}algebra iff $\jj$ is a simple Lie {\it super}algebra.

The {\it derived algebra} $\rr'$ of a conformal superalgebra
$\rr$ is defined as the span over $\ccc$ of all elements
$a_{(j)}b$, where $a,b\in\rr$, $j\in\zz_+$.  It is easy to see
that this is an ideal of $\rr$ such $\rr/\rr'$ is a commutative
conformal superalgebra (i.e. all products of $\rr/\rr'$ are
zero).  One defines $\rr''=(\rr')'$, etc, obtaining a descending
sequence of ideals $\rr\supset\rr'\supset\rr''\supset \ldots$. The
conformal superalgebra $\rr$ is called {\it solvable} if $n$-th
member of this sequence is zero for some $n>0$. The conformal
superalgebra $\rr$ is called {\it semisimple} if its only
solvable ideal is zero.

Define by induction another descending sequence of ideals
$\rr\supset\rr^{1} \supset\rr^{2}\supset \ldots$ by letting
$\rr^{1}=\rr', \ldots, \rr^{n}=$ span over $\ccc$ of all products
$a_{(j)}b$ where $a\in\rr$, $b\in\rr^{n-1}$, $j\in\zz_+$.  The
conformal superalgebra $\rr$ is called {\it nilpotent} if $n$-th
member of this sequence is zero for some $n>0$.

\section{Preliminaries on conformal modules.}
\label{sec:2}

Let $(\jj,\ff)$ be a Lie superalgebra of formal distributions, and let $V$ 
be a $\jj$-module. We say that a formal distribution $a(z)\in\ff$ and a 
formal distribution $v(z)=\sum_{n\in\zz}v_{(n)}z^{-n-1}$ with coefficients 
in $V$ are {\it local} if 
\begin{equation}
  (z-w)^Na(z)v(w)=0 \text{ for some } N\in\zz_+ \, .
  \label{2.1}
\end{equation} 
In the same way as in \cite{1}, \S{2.3}, one shows that \req{2.1}
is equivalent to 
\begin{equation} 
  a(z)v(w)=\sum_{j=0}^{N-1}(a_{(j)}v)(w)\pa_w^j\de(z-w)/j! \, ,
  \label{2.2}
\end{equation}
for some formal distributions $(a_{(j)}v)(w)$ with coefficients in $V$, which
are uniquely determined by the formula
\begin{equation}
  (a_{(j)}v)(w)=\Res_z(z-w)^ja(z)v(w) \, .
  \label{2.3}
\end{equation}

%%%\examp{Example 2.1} 
\begin{example}
\label{ex:2.1}
Consider the following representation of the (centerless)
Virasoro algebra in the vector space $V$ with the basis $v_{(n)}$, $n\in\zz$,
over $\ccc$:
\begin{equation*}
  L_mv_{(n)}= \left( (\Delta-1) ( m+ 1) -n \right)
  v_{(m+n)} + \alpha v_{(m+n+1)} \, ,
\end{equation*}
where $\Delta,\alpha\in\ccc$. In terms of formal distributions $L(z)$ and 
$v(z)=\sum_nv_{(n)}z^{-n-1}$ this can be written as follows:
\begin{equation}
  L(z)v(w)=(\pa+\alpha)v(w)\de(z-w)+\Delta v(w)\de_{w}'(z-w) \, .
  \label{2.4}
\end{equation}
Hence $L(z)$ and $v(z)$ are local.
\end{example}

Suppose that $V$ is spanned over $\ccc$ by coefficients of a family $E$ of
formal distributions such that all $a(z)\in\ff$ are local with all 
$v(z)\in E$.
Then we call $(V,E)$ a {\it conformal module over} $(\jj,F)$.

The following is a representation-theoretic analogue (and a generalization)
of Dong's lemma (see \cite{1},~\S{3.2}).

\begin{lemma} 
  \label{lem:2.1}
Let $V$ be a module over a Lie superalgebra $\jj$, let $a(z)$
  and $b(z)$ (resp.\ $v(z)$) be formal distributions with
  coefficients in $\jj$ (resp.~$V$). Suppose that all the pairs
  $(a,b),(a,v)$ and $(b,v)$ are local.  Then the pairs
  $(a_{(j)}b,v)$ and $(a,b_{(j)}v)$ are local for all
  $j\in\zz_+$.
\end{lemma}

\pr{Proof}
%\noindent{\it Proof}. 
We may assume that all three pairs satisfy respectively (\ref{1.1}) 
and \req{2.1} for some $N\in\zz_+$. Then we have:
\begin{equation*}
\begin{split}
&(z-w)^{3N}(a_{(j)}b)(z)v(w) \\
&\qquad {} = (z-w)^N \Res_u
\sum_{i=0}^{2N}\binom{2N}{i}(z-u)^i(u-w)^{2N-i}
(u-z)^j[a(u),b(z)]v(w) \, .
\end{split}
\end{equation*}
The summation over $i$ in the right-hand side from $0$ to $2N$ may be 
replaced by that from $0$ to $N$ since $a(u)$ and $b(z)$ are mutually local.
Hence it can be written as follows:
\begin{equation*}
(z-w)^N\Res_u(u-w)^N P(z,u,w) (u-z)^j  \left(a(u) b(z) v(w) -
 b(z) a(u) v(w) \right)
\end{equation*} 
for some polynomial $P$. But this is zero since both pairs $(b,v)$
and $(a,v)$ are local, which proves that the pair $(b_{(j)}a,v)$ is
local.

Next, using the first part of lemma, we may find $N$ for which all pairs 
$(a_{(j)}b,v)$ and $(a,v)$ satisfy \req{2.1}.
Then we have:
\begin{equation*}
\begin{split}
a(z)(b_{(j)}v)(w) &= \Res_ua(z)b(u)v(w)(u-w)^j\\
                 &=-\Res_u([b(u),a(z)]v(w)+b(u)a(z)v(w))(u-w)^j\\
                 &=-\Res_u \left( \sum_{i\geq0}
                 (b_{(i)}a)(z)v(w)\pa_z^i\de(u-z)/i! + 
                 b(u)a(z)v(w) \right)
                 (u-w)^j     \, .
\end{split}
\end{equation*} 
Hence $(z-w)^Na(z)(b_{(j)}v)(w)=0$.  \Qed

This lemma shows that the family $E$ of a conformal module $(V,E)$ over 
$(\jj,\ff)$ can always be included in a larger family $E^c$ which is
local with respect to $\ff^c$ and such that $\pa E^c\subset E^c$ and 
$a_{(j)}E^c\subset E^c$
for all $a\in\ff$ and $j\in\zz_+$.

It is straightforward to check the following properties for $a,b\in\ff$
and $v\in E^c$:
\begin{eqnarray}  
  \left[a_{(m)},b_{(n)}\right]v
  &=&
  \sum_{j=0}^m\binom{m}{j}\left(a_{(j)}b\right)_{(m+n-j)}v,
  \label{2.5} \\
  (\pa a)_{(n)}v &=& \left[\pa,a_{(n)}\right]v=-na_{(n-1)}v.
  \label{2.6}
\end{eqnarray}
(Here [~,~] is the bracket of operators on $E^c$.) It follows
from \req{2.5} (by  
induction on $m$) and \req{2.6}, that $a_{(j)}E^c\subset E^c$ for all
$a \in F^c$
and $j\in\zz_+$.  

Thus, any conformal module $(V,E)$ over a Lie superalgebra of formal 
distributions $(\jj,\ff)$ gives rise to a {\it module} $M=E^c$ over the 
conformal algebra $R=\ff^c$, defined as follows. It is a (left) $\zz/2\zz$-
graded $\ccc[\pa]$-module with $\ccc$-linear maps $a\longmapsto a_{(n)}^M$
of $R$ to $End_\ccc M$ given for each $n\in\zz_+$ such that the following 
properties hold (cf.~\req{2.5} and \req{2.6}) for $a,b\in\rr$,
$m$, $n\in\zz_+$: 
%\begin{equation*}
%%%\begin{gather*}
\begin{deflist}{M}
\setcounter{bean}{-1}
\item
  $a_{(n)}^Mv=0 \text{ for } v\in M \text{ and } n\gg 0$,
  %%%\tag{M0}\\
\item
  $\left[a_{(m)}^M,b_{(n)}^M\right]=
  \sum_{j=0}^m\binom{m}{j}\left(a_{(j)}b\right)_{(m+n-j)}^M$,
  %%%\tag{M1}\\  
\item
  $(\pa a)_{(n)}^M=\left[\pa,a_{(n)}^M\right]=-na_{(n-1)}^M$.
  %%%\tag{M2}   
\end{deflist}
%%%\end{gather*}
%\end{equation*}

Conversely, suppose that a conformal superalgebra $\rr$ is a free 
$\ccc[\pa]$-module and consider the associated Lie superalgebra of
formal distributions $(\jj(\rr),\ff)$ (see \S\ref{sec:1}). Let
$M$ be a module over  
the conformal superalgebra $\rr$ and suppose that $M$ is a free 
$\ccc[\pa]$-module with basis $\{v^{\alpha}\}_{\alpha\in J}$. This gives rise
to a conformal $\jj(\rr)$-module $V(M)$ with basis $v_{(n)}^{\alpha}$, where
$\alpha\in J$, $n\in\zz$, defined by (cf.~\req{2.2}):
\begin{equation}
  a^{\alpha}(z)v^{\beta}(w)=\sum_{j\in\zz_+}\left(a_{(j)}^{\alpha}v^{\beta}
  \right)(w)\pa_w^j\de(z-w)/j!.
  \label{2.7}
\end{equation} 
%
%%%\noindent{\bf Remark 2.1.} 
\begin{remark}
\label{rem:2.1}
Given a module $M$ over a conformal superalgebra $\rr$, we may change its 
structure as a $\ccc[\pa]$-module replacing $\pa$ by $\pa+A$ where $A$ is an 
endomorphism over $\ccc$ of $M$ which commutes with all $a_{(n)}^M$ (this
will not affect axiom (M2)). Sometimes we shall not distinguish these
$\rr$-modules. For example, we may put $\alpha=0$ in~\req{2.4}.
\end{remark}

A conformal module $(V,E)$ (resp.\ module $M$) over a Lie superalgebra of 
formal distributions $(\jj,\ff)$ (resp.\ over a conformal superalgebra $\rr$)
is called {\it finite} if $E^c$ (resp.~$M$) is a finite $\ccc[\pa]$-module. 

Note that the maps a $a\mapsto a_{(n)}$ of $\rr$ to $End_{\ccc}\rr$
define a $\rr$-module, called the {\it adjoint module}. (It is finite iff 
$\rr$ is finite).

The following lemma generalizes Proposition 2.7 from \cite{1}.
\begin{lemma}
\label{lem:2.2} 
Let $M$ be a $\rr$-module.
\vspace*{-1.3ex}
\alphaparenlist
\begin{enumerate}
\item %%%{\rm(a)} 
If $\pa v= \lambda v$ for some $\lambda\in\ccc$ and 
$v\in M$,
then $\rr_{(n)}^Mv=0$ for all $n\in\zz_+$.
\vspace*{-1.3ex}
\item %%%{\rm(b)} 
If $M$ is a finite module which has no non-zero invariants, i.e. 
vectors
$v$ such that $\rr_{(n)}^Mv=0$ for all $n\in\zz_+$, then it is a free 
$\ccc[\pa]$-module.
\end{enumerate}
\end{lemma}

\pr{Proof} Let $a\in\rr$ and take the minimal $m$ such that 
$a_{(j)}v=0$ for $j\geq m$.
We have: $0=(\pa-\lambda)\left(a_{(m)}v\right)=
\left[\pa-\lambda,a_{(m)}\right]v+a_{(m)}
(\pa-\lambda)v=-ma_{(m-1)}v$. This shows that $m=0$, proving (a).

Since (a) says that a $\rr$-module $M$ with only zero invariants has no 
torsion, (b) follows from (a). 
\Qed

The correspondence between finite conformal modules with no non-trivial 
invariants over a finite Lie superalgebra of formal distributions 
$(\jj,\ff)$ with a trivial center and finite conformal modules with no 
non-trivial invariants over the conformal superalgebra $\ff^c$ is described
in a fashion similar to Proposition~1.1.

Note that Example 2.1 gives a 2-parameter family of (irreducible) modules
of size (1,0) over the Virasoro conformal algebra. Note also that the 
well-known family of graded $\Vir$-modules given by
\begin{equation*}
L_mv_{(n)}=((\Delta-1)m-n+\alpha)v_{(m+n)}
\end{equation*}
is conformal, but is finite iff $\alpha = \Delta - 1$.

The following simple observation is fundamental for representation theory
of conformal superalgebras.

\begin{proposition}
\label{prop:2.1}
Consider the Lie superalgebra of formal distributions $(\jj(\rr),\ff)$
defined by (\ref{1.5}) and let $\jj(\rr)_+$ be the subalgebra of $\jj(\rr)$
spanned by all $a_{(n)}^i$ with
$i\in I$, $n\in\zz_+$. Denote
by $\jj(\rr)^+$
the semidirect product of the $1$-dimensional Lie
algebra $\ccc\pa$ and
ideal $\jj(\rr)_+$ with the action of $\pa$ on $\jj(\rr)_+$ given
by 
$\pa(a_{(n)}^i)=-na_{(n-1)}^i$. Then any module $M$ over the conformal
superalgebra $\rr$ gives rise to a $\jj(\rr)^+$-module $M$ (over $\ccc$) 
such that 
\begin{equation}
  a_{(n)}^i v=0 \text{ for } v\in M \text{ and } n\gg 0 \, .
  \label{2.8}
\end{equation}
\end{proposition}

\begin{cor}
\label{cor:2.1} 
Let $\rr\!=\!\oplus_{\alpha\in I}\ccc[\pa]a^\alpha$ be a conformal 
superalgebra and 
$M\!=\!\oplus_{\beta\in J}\ccc[\pa]v^{\beta}$ be a free
$\ccc[\pa]$-module. Then,
given $a_{(n)}^{\alpha}v^{\beta}\in M$ for all $\alpha\in I$, $\beta\in J$,
$n\in\zz_+$, which is $0$ for $n\gg 0$, we may extend uniquely the action 
of $a_{(n)}^{\alpha}$ to the whole $\rr$ on $M$ using {(M2)}.
Suppose that {(M1)} holds for all $a=a_{(m)}^{\alpha}$,
$b=a_{(n)}^{\beta}$. Then $M$ is a $\rr$-module.
\end{cor}

Using Proposition~\ref{prop:2.1} and Corollary~\ref{cor:2.1}, one can
construct large families of finite modules over conformal
superalgebras.

%\smallskip
%%%\examp{Example 2.2} 
\begin{example}
\label{ex:2.2}
Let $\Vir_{\geqslant 0}=\sum_{j\geqslant 0}\ccc L_j$ and
consider a representation $\pi$ of  $\Vir_{\geqslant 0}$ in a 
finite-dimensional (over $\ccc$) vector space $U$.
Let $A$ be an endomorphism of $U$ commuting with all 
$\pi(L_j)$ $(j\in\zz_+)$. Then $\ccc[\pa]\otimes U$ is a finite module 
over the conformal algebra $\Vir$ defined by the following formulas ($u\in U$):
\begin{equation*}
L_{(0)}u=(\pa+A)u \, , \quad 
L_{(j)}u=\pi(L_{j-1})u \text{ for } j\geqslant 1 \, .
\end{equation*}
For example, we can take $\pi(L_0)=B$, where $B$ is an endomorphism of $U$
commuting with $A$. Then 
\begin{equation*}
L_{(0)}u=(\partial+A)u \, , \quad 
L_{(1)}u=Bu \, , \quad 
L_{(j)}u=0 \, , \text{ for } j>1 \, ,
\end{equation*}
defines a finite module over $\Vir$, which we denote by $M(A,B)$.
Taking $\ddim U=1$, $A=\alpha$ and $B=\Delta$ gives Example~\ref{ex:2.1}.
\end{example}

\section{Structure theory of finite conformal algebras.}
\label{sec:3}

Results stated in this section is a joint work \cite{3} with my student, 
Allessandro\break 
D'Andrea. The reader is referred to \cite{3} for
details.

\begin{theorem} {\rm(conformal Engel theorem)} 
\label{theor:3.1}
Let $M$ be a finite module over
a finite conformal (super)algebra $\rr$ such that $a_{(n)}^M$ is a
nilpotent operator for all $a\in\rr$ and $n\in\zz_+$. Then there exists
a non-zero vector $v\in M$ such that
\begin{equation}
  a_{(n)}^Mv=0 \text{ for all } a\in\rr \, , \; n\in\zz_+ \, .
  \label{3.1}
\end{equation} 
\end{theorem}

\begin{cor}
\label{cor:3.1}
A finite conformal (super)algebra $\rr$ is nilpotent iff $a_{(n)}$ is a 
nilpotent operator (on $\rr$) for all $a\in\rr$ and $n\in\zz_+$.
\end{cor}

\begin{theorem}
{\rm (conformal Lie theorem)} 
\label{theor:3.2}
Let $M$ be a finite module over a finite 
solvable 
conformal algebra $\rr$. Then there exists a non-zero vector $v\in M$
such that 
\begin{equation}
  a_{(n)}^Mv=\lambda(a,n)v \, , 
  \text{ where } \lambda(a,n)\in\ccc \, , 
  \text{ for all } a\in\rr \, , \; n\in\zz_+ \, . 
  \label{3.2}
\end{equation}
\end{theorem}

\begin{cor}
\label{cor:3.2}
The derived algebra of a finite solvable conformal algebra is nilpotent.
\end{cor}

\begin{theorem}
\label{theor:3.3}
A simple finite conformal algebra is isomorphic either to a
current conformal algebra $\rr(\jjb)$, where $\jj$ is a simple
finite-dimensional Lie algebra, or to the Virasoro conformal
algebra $\Vir$.
\end{theorem}

It is not quite true that a semisimple finite conformal algebra is a direct 
sum of simple ones. For example, if $\jj$ is a semisimple finite-dimensional
Lie algebra, the semi-direct sum $\rr=\Vir + \rr (\jjb)$ is semisimple,
where the products are given by:
\begin{equation*}
\begin{split}
L_{(0)}a &= \pa a \text{ if } a \in \rr ; \quad 
            L_{(1)} L = 2L;  \quad
            L_{(1)} a = a \text{ if } a \in \jj ; \\ 
a_{(0)}b &= [a,b]; \quad \text{all other products (up to the
  order) on } L + \jj \text{ are } 0 .
\end{split}
\end{equation*} 

\begin{theorem}
\label{theor:3.4}
A semi-simple finite conformal algebra is isomorphic to the direct sum of
conformal algebras of three types:
\vspace*{-1.3ex}
\alphalist
\begin{enumerate}
\item %%%{\rm a)} 
current conformal algebra $\rr(\jjb)$,
where $\jj$ is simple, 
\vspace*{-1.3ex}

\item %%%{\rm b)} 
conformal algebra $\Vir$, 
\vspace*{-1.3ex}

\item %%%{\rm c)} 
conformal algebra 
$\Vir + \rr(\jjb)$, where $\jj$ is semisimple.
\end{enumerate}
\end{theorem}

\begin{theorem}
\label{theor:3.5}
A finite conformal algebra $\rr$ admits a faithful finite
irreducible module iff $\rr$ is of one of the following types:
\vspace*{-1.3ex}
\alphalist
\begin{enumerate}
\item %%%{\rm a)} 
current conformal algebra $\rr(\jjb)$, where $\jj$ is semisimple,
\vspace*{-1.3ex}

\item %%%{\rm b)} 
current conformal algebra $\rr(\jjb)$, where $\jj$ is semisimple or zero
plus $1$-dimen\-sional,
\vspace*{-1.3ex}

\item %%%{\rm c)} 
conformal algebra $\Vir + \rr(\jjb)$, where $\jj$ is semisimple or zero,
\vspace*{-1.3ex}

\item %%%{\rm d)} 
conformal algebra $\Vir + \rr(\jjb)$, where $\jj$ is semisimple or zero 
plus $1$-dimensional,\vspace*{-1.3ex}

\item %%%{\rm e)} 
one of {a)} or {c)} plus $1$-dimensional ($\text{over } \ccc$) 
conformal algebra. 
\end{enumerate}
\end{theorem}

\section{Finite modules over semisimple finite conformal algebras.}
\label{sec:4}

Results of this section is a joint work \cite{4} with Sun-Jen Cheng.

The following is a key lemma.

\begin{lemma}
  \label{lem:4.1}
  Let ${\cal L}$ be a Lie superalgebra (over $\ccc$) with a
  distinguished element $\pa$ and a descending sequence of
  subspaces ${\cal L}\supset{\cal L}_0 \supset{\cal L}_1\supset
  \ldots$ such that $[\pa,{\cal L}_n]={\cal L}_{n-1}$ for all
  $n>0$. Let $V$ be a ${\cal L}$-module and let
\begin{equation*}
V_n=\{v\in V\mid {\cal L}_nv=0\}, n\in\zz_+ \, .
\end{equation*}                                                         
Suppose that $V_n\neq0 \text{ for } n\gg0$.  Suppose that the
minimal $N\in\zz_+$ for which $V_N\neq0$ is positive.  Let
$v_1,v_2,\ldots$ be a linearly independent over $\ccc$ set of
vectors of $V_N$ which generate $\ccc[\pa]V_N$ as a
$\ccc[\pa]$-module. Then $v_1,v_2,\ldots$ is a basis 
$(\text{over } \ccc)$ of $V_N$ and a free set of generators of
the $\ccc [\pa]$-module
$\ccc[\pa]V_N$. In particular, $V_N$ is finite-dimensional over
$\ccc$ if V is a finitely generated $\ccc[\pa]$-module.
\end{lemma}

It is clear from the very definition that the Lie superalgebra
${\cal L}=\jj(\rr)^+$ and any finite conformal module $V$ over
$\rr$ viewed as an ${\cal L}$-module (see
Proposition~\ref{prop:1.1}) satisfy the conditions of
Lemma~\ref{lem:4.1}.

\begin{theorem}
\label{theor:4.1}
Any irreducible finite $\Vir$-module is either $1$-dimensional
over $\ccc$, or else is a free of rank $1$ $\ccc[\pa]$-module 
$\ccc[\pa]v$  defined by (cf.\ Example 2.1):
\begin{equation*}
L_{(0)}v=(\pa+\alpha)v \, , \quad L_{(1)}v=\Delta v \, , \quad
L_{(j)}v=0 \text{ for } j>1 \, 
\end{equation*}  
for some $\alpha$, $\Delta\in\ccc$, $\Delta\neq0$.
\end{theorem} 

Let $\jj$ be a semisimple finite-dimensional Lie algebra and let $\rr(\jjb)$
be the associated current conformal algebra.

\begin{theorem}
\label{theor:4.2}
Any irreducible finite $\rr(\jjb)$-module is either $1$-dimensional over 
$\ccc$,
or else is of the form $\ccc[\pa]\otimes_{\ccc}U$,where $U$ is a non-trivial 
irreducible finite-dimensional $\jj$-module, with the action of $\rr(\jjb)$
defined by:
\begin{equation*}
a_{(j)}u=0 \text{ for } j>0 \, , \quad
a_{(0)}u=au \, , \text{ if } a\in\jj \, , \; u\in U.
\end{equation*}
\end{theorem}

\begin{theorem}
\label{theor:4.3}
  Any irreducible finite $\Vir + \rr(\jjb)$-module is either the
  trivial $1$-dimen\-sional over $\ccc$, or else is of the
  form $\ccc[\pa]\otimes_{\ccc} U$, where $U$ is a non-trivial
  irreducible finite-dimensional $\jj$-module, with the action of
  $\Vir + \rr(\jjb)$ defined by:
\begin{equation*}
\begin{split}
&a_{(j)}u=\delta_{j0}au \text{ if } a\in\jj \, ,  \; u\in U \, ,
\; j\in\zz_+ \, , \\
&L_{(0)}u=\pa u \, , \quad
L_{(j)}u= \delta_{1j}\Delta u \,  \text{ if } u\in U \, , \; 
j\geq1 \, , \text{ where } \Delta\in\ccc \, .  
\end{split}
\end{equation*}
\end{theorem}

Thus, we see that the classification of irreducible finite modules over 
semi-simple conformal algebras is similar to that of semi-simple Lie 
algebra, a small difference being caused by the Virasoro conformal 
algebra. 

There is, however, an essentially new feature: 
complete reducibility breaks down 
in a very dramatic way. For example, the module $M(A,B)$ over $\Vir$ 
constructed in Example~2.2 is indecomposable iff the pair of commuting
operators $(A,B)$ is indecomposable. Similar examples may be constructed
for current conformal algebras as well.

We describe below all the extensions between an irreducible $\Vir$-module
$M(\alpha,\Delta)$, $\alpha,\Delta\in\ccc, \Delta\neq 0$, and the
trivial $\Vir$-module $\ccc$.

\begin{theorem}
\label{theor:4.4}
{\rm(a)} A non-split extension of $\Vir$-modules of the form
\begin{equation*}
0\longrightarrow M(\alpha,\Delta)\longrightarrow V\longrightarrow\ccc
\longrightarrow 0
\end{equation*}
$\text{exists iff } \Delta=1 \, \text{ which is }$
\begin{equation*}
0\longrightarrow M(\alpha,1)\longrightarrow M(\alpha,0)\longrightarrow\ccc
\longrightarrow 0 \, .
\end{equation*}
{\rm(b)} A non-split extension of $\Vir$-modules of the form
\begin{equation*}
0\longrightarrow\ccc\longrightarrow V\longrightarrow M(\alpha,\Delta)
\longrightarrow 0
\end{equation*}
exists iff $\Delta=1$ or $2$. In these two cases $V=\ccc[\pa]v+\ccc$, where 
$\ccc$ is the trivial submodule, with the action on $v$:
\begin{equation*}
L_{(0)}v=(\pa+\alpha)v \, , \quad 
L_{(1)}v=\Delta v \, , \quad
L_{(j)}v=\delta_{j,\Delta+1}1\in \ccc \text{ for } j>1 \, .
\end{equation*} 
{\rm(c)} \cite{18} A non-split extension of $\Vir$-modules of the form
\begin{equation*}
0\longrightarrow M(\alpha,\Delta)\longrightarrow V\longrightarrow 
M(\beta,\Delta')\longrightarrow 0
\end{equation*}
may exist only when $\alpha=\beta$ and $\Delta'-\Delta=0,1,2,3,4,5,6$. They
can be explicitly described using the results of~\cite{5}.
\end{theorem}

Denote by $M(\lambda)$ the $\rr(\jjb)$-module $\ccc[\pa]\otimes F(\lambda)$,
where $F(\lambda)$ is the finite-dimen\-sional irreducible highest weight
module over the simple Lie algebra $\jj$ (see Theorem~\ref{theor:4.2}).  
                      
\begin{theorem}
\label{theor:4.5}
The only non-split extension between $1$-dimensional $\rr(\jjb)$-module
and a module $M(\lambda)$ is of the form:
\begin{equation*}
0\longrightarrow\ccc\longrightarrow V\longrightarrow M(\theta)
\longrightarrow 0,  
\end{equation*}
where $\theta$ is the highest root (so that $F_{\theta}=\jj$). In this case 
$V=\ccc[\pa]\otimes\jj+\ccc$, where $\ccc$ is a submodule, with the action 
on $\jj$:
\begin{equation*}
a_{(0)}b=[a,b] \, , \quad
a_{(j)}b=\delta_{1j}(a\mid b)\in\ccc \text{ for } j\geq 1 \, ,
\end{equation*} 
where $(.|.)$ is the Killing form on $\jj$.
\end{theorem}

\section{Classification of simple finite conformal superalgebras.}
\label{sec:5}

The list of all finite conformal superalgebras is much richer
than that of finite conformal algebras. First, there are many
more simple finite-dimensional Lie superalgebras (classified in
\cite{6}), and the associated conformal superalgebra is finite and
simple. Second, there are many ``superizations'' of the Virasoro
conformal algebra. We describe them below.  They are associated
to superconformal algebras constructed in \cite{6}, \cite{7} and
\cite{8} (cf.~\cite{9}).

Let $\Lambda(N)$ denote the Grassmann algebra over $\ccc$ in $N$ 
indeterminates $\xi_1,\ldots,\xi_N$. Let $W(N)$ be the Lie superalgebra of 
all derivations of superalgebra $\Lambda(N)$. It consists of all linear 
differential operators $\sum_{i=1}^{N}P_i\pa_i$, where $P_i\in\Lambda(N)$
and $\pa_i$ stands for the partial derivative by $\xi_i$.

The first series of examples is the series conformal superalgebras $W_N$
of rank $(N+1)2^N$:
\begin{equation*}
W_N=\ccc[\pa]\otimes_{\ccc}(W(N)\oplus\Lambda(N))
\end{equation*}    
with the following products $(a,b\in W(N), \, f, \, g\in\Lambda(N))$:
\begin{equation*}
\begin{split}
&a_{(j)}b=\delta_{j0}[a,b] \, , \quad
a_{(0)}f=a(f) \, , \quad
a_{(j)}f=-\delta_{j1}(-1)^{p(a)p(f)}fa \text{ if } j\geq 1 \, , \\
&f_{(0)}g=-\pa(fg) \, ,\quad
f_{(j)}g=-2\delta _{j1}fg \text{ if } j\geq 1 \, .
\end{split}
\end{equation*}

For an element $D=\sum_{i=1}^NP_i(\pa,\xi)\pa_i+f(\pa,\xi)\in W_N$ define
{\it divergence} by the formula
\begin{equation*}
\ddiv D=\sum_{i=1}^N(-1)^{p(P_i)}\pa_iP_i+\pa f \, .
\end{equation*}

The second series of examples is 
\begin{equation*}
S_N=\{D\in W_N\mid \ddiv D=0\}.
\end{equation*}
This is a conformal superalgebra of rank $N2^N$, and it is simple iff $N\geq2$.

The third series of examples is $K_N$. It is also a subalgebra of $W_N$ 
(of rank $2^N$), but it is more convenient to describe it as follows:
\begin{equation*}
K_N=\ccc[\pa]\otimes_\ccc\Lambda(N)
\end{equation*}
with the following products $(f,g\in\Lambda(N))$:
\begin{equation*}
\begin{split}
&f_{(0)}g=\left(\frac{1}{2}|f|-1\right)\pa fg+\frac{1}{2}(-1)^{|f|}\sum_{i=1}
^N(\pa_if)(\pa_ig) \, ,\\
&f_{(j)}g=\left(\frac{1}{2}(|f|+|g|)-2\right)\de_{j1}fg 
\text{ if } j\geq1 \, .
\end{split}
\end{equation*}
We assume here that $f$ and $g$ are homogeneous elements of degrees $|f|$
and $|g|$ in the gradation defined by $\ddeg\xi_i=1$ for all $i$.

These three series include all well-known examples. Thus, $W_0\backsimeq K_0$
is the Virasoro conformal algebra, $K_1$ is the Neveu-Schwarz conformal 
superalgebra, $K_2\backsimeq W_1$ and $K_3$ are the $N=2$ and $3$ conformal 
superalgebras, $S_2$ is the $N=4$ conformal superalgebra. Nevertheless
there is one exceptional example constructed in \cite{8}. It is the subalgebra
$CK_6$ of rank $32$ in the conformal superalgebra $K_6$ spanned over 
$\ccc[\pa]$ by the following elements:
\begin{equation*}
1+\alpha\pa^3\nu \, , \quad 
\xi_i-\alpha\pa^2\xi_i^* \, , \quad
\xi_i\xi_j+\alpha\pa(\xi_i\xi_j)^* \, ,  \quad
\xi_i\xi_j\xi_k+\alpha(\xi_i\xi_j\xi_k)^* \,.
\end{equation*}
Here $\alpha^2=-1$, 
$\nu=\xi_1\xi_2\ldots\xi_6$,
$(\xi_{i_1}\xi_{i_2}\ldots)^*= \pa_{i_1}\pa_{i_2}\ldots \nu$.

%%%\rema{Remark 5.1} 
\begin{remark}
\label{rem:5.1}
Let $\jj_1'=\jj(CK_6)$ (resp.\ $\jj=\jj(K_6)$)
be the Lie superalgebra associated to the conformal superalgebra
$CK_6$ (resp.\ $K_6$) and let $\jj'_+$ (resp.\ $\jj_+$) be its
subalgebra spanned by all non-negative modes of all formal
distributions of $CK_6$ (resp.\ $K_6$). The Lie superalgebra
$\jj_+$ has a $\zz$-gradation of the form
$\oplus_{j\geq{-2}}\jj_j$ with $\jj_{-2}=
\ccc$, $\jj_{-1}=\ccc^6$, $\jj_0=so_6\oplus\ccc$,
$\jj_1=\bigwedge^3\ccc^6\oplus\ccc^6$ \cite{6}.  The reason for
occurrence of $CK_6$ lies in the fact that the $so_6$-module
$\Lambda^3\ccc^6\backsimeq V\oplus V^*$ is reducible, where $V$
is a $10$-dimensional $so_6$-module. The subspace $\jj_{-1}\oplus
\jj_0\oplus (V+\ccc^6)$ generates the simple $\zz$-graded
superalgebra $\jj_+'$, which should be added to the list of
Theorem~10 from~\cite{6}.
\end{remark}

\begin{theorem}
\label{theor:5.1}
Any simple finite conformal superalgebra is isomorphic either to the current 
conformal superalgebra associated to a simple finite-dimensional Lie
superalgebra (classified in \cite{6}), or to one of the following conformal
superalgebras $(N\in\zz_+)$:
\begin{equation*}
W_N \, , \quad 
S_{N+2}\, , \quad
K_N \, , \quad 
CK_6 \, .
\end{equation*}
\end{theorem}

The proof of this theorem is pretty long \cite{10}. It follows the same 
strategy as \cite{6}, the first step have been worked out in \cite{9}.

\section{Work in progress.}
\label{sec:6}

%6.1
\subsection{Cohomology \protect\cite{11}} 
\label{sec:6.1}

In representation theory and conformal field 
theory it is important to consider central extensions of Lie (super)algebras
of formal distributions. In the language of a conformal (super)algebra $\rr$
it amounts to considering a conformal superalgebra
\begin{equation*}
\rrt=\rr\oplus\ccc C \text{ with } p(C)=0 \, , \quad
\pa C=0 \, , \quad
C_{(n)}\rrt=0
\text{ for } n\in\zz_+ \, .
\end{equation*}
The $n$-th product $a_{(\nt)}b$ on $\rr\subset\rrt$ is given by
\begin{equation*}
a_{(\nt)}b=a_{(n)}b+\alpha_n(a,b)C \, .
\end{equation*}
It is straightforward to see that $\rrt$ is a conformal superalgebra 
iff the sequence $\{\alpha_n\}_{n\in\zz_+}$ is a $2$-cocycle on $\rr$ 
defined as follows. It is a sequence of $\ccc$-valued $\ccc$-bilinear 
forms $\alpha_n(n\in\zz_+)$ on $\rr\times\rr$ such that $(a,b,c\in\rr, m,n
\in\zz_+)$:
\begin{equation*}
\begin{split}
\alpha_n(\pa a,b) &= -n\alpha_{n-1}(a,b) \, ,\\
\alpha_n(a,b) &= (-1)^{n+1+p(a)p(b)}\alpha_n(b,a) \, ,\\
\alpha_m\left(a,b_{(n)}c\right) &=
\sum_{j=0}^m\binom{m}{j}\alpha_{m+n-j}\left(a_{(j)}b,c\right)+(-1)
^{p(a)p(b)}\alpha_n\left(b,a_{(m)}c\right) \, .
\end{split}
\end{equation*}

As usual, the trivial cocycle $\alpha_n(a,b)=f\left(a_{(n)}b\right)$, 
where $f:\rr
\rightarrow\ccc$ is a $\ccc$-linear map, defines a trivial central extension
of $\jj$ (isomorphic to the direct sum of $\jj$ and~$\ccc$). Two cocycles that
differ by a trivial cocycle are called equivalent; they define isomorphic 
central extensions.

Starting with this observation, Alexander Voronov and I developed
in \cite{11} a cohomology theory 
of conformal algebras. This cohomology $H^*(\rr,M)$ is defined for any 
$\rr$-module $M$ in a way similar to the Lie algebra cohomology. The central
extensions of $\rr$ are parameterized by $H^2(\rr,\ccc)$. 

%%%\goodbreak\vbox{
%%%\examp{Examples}
\goodbreak
\begin{examples}
\vspace*{-1.3ex}
\alphalist
\begin{enumerate}
\item %%%(a) 
$\ddim_{\ccc}H^2(\rr(\jjb),\ccc)=1$ if $\jj$ is a simple finite-dimensional
Lie algebra. The non-trivial $2$-cocycle is given by
\begin{equation*}
\alpha_n(a,b)=\delta_{n,1}(a\mid b) \text{ for } a,b\in\jj \, .
\end{equation*}\vspace*{-4ex}

\item %%%(b) 
\cite{7}, \cite{8}. Let $\rr$ be one of the simple conformal 
superalgebras $W_N, S_N, K_N, CK_6$. Then $\ddim_{\ccc} H^2(\rr,\ccc)=0$ or 
$1$ and the latter case takes place iff $\rr=W_N$ with $N\leq2$, $S_2$, and 
$K_N$ with $N\leq4$.
\vspace*{-1.3ex}

\item %%%(c) 
\cite{11} $\ddim_{\ccc}H^n(\Vir,\ccc)=1$ if $n=0,2$ or $3$, and $=0$ 
otherwise. This example is intimately related to the Gelfand-Fuchs 
cohomology \cite{5}.    
\end{enumerate}
%%%}
\end{examples}

%6.2
\subsection{Modules over finite simple conformal superalgebras 
\protect\cite{12}}
\label{sec:6.2}

Combining methods described in \S\ref{sec:4} with the ones of
\cite{13} Alexey Rudakov and I were able to classify and
construct all finite irreducible modules over all finite simple
conformal superalgebras. For current conformal superalgebras and
for $K_1$ (Neveu-Schwarz) it was done already in \cite{4}, and
the answer is similar to Theorem \ref{theor:4.1} and
\ref{theor:4.2} respectively. However, in all other cases some
interesting new effects occur.

%6.3
\subsection{Finite modules over infinite conformal algebras
  \protect\cite{15}}
\label{sec:6.3}

The next after finite conformal (super)algebras is the class of conformal
algebras $\rr$ that admit a faithful finite module $M$. Suppose for 
simplicity that $M$ is a free module of rank $N$ over $\ccc[\pa]$. It 
turns out  that $\rr$ is then a subalgebra of the ``general''
conformal algebra $\gcn$ of all ``conformal'' endomorphisms of $M$. A 
conformal endomorphism $A$ of $M$ is a sequence of endomorphisms $A_{(n)}$
of $M$ over $\ccc$ for each $n\in\zz_+$ such that $A_{(n)}v=0$ for
$n \gg 0$ and 
%\begin{equation*}
$A_{(n)}\pa v=\pa A_{(n)}v+nA_{(n-1)}v, v\in M
$.
%\end{equation*}                                                   
All conformal endomorphisms of $M$ form an infinite conformal algebra, 
denoted by $\gcn$, with $\ccc [\partial]$-module structure 
$(\pa A)_{(n)}=-nA_{(n-1)}$ and the following products for each 
$m\in\zz_+$:
\begin{equation*}
(A_{(m)}B)_{(n)}=\sum_{j=0}^m(-1)^{m+j}\binom{m}{j}[A_{(j)},B_{(m+n-j)}],
\end{equation*} 
$M$ being its finite irreducible module. This conformal algebra is simple, 
but infinite.

It is easy to see that $\gcn$ is the conformal algebra associated to the Lie
algebra of 
$N\times N$ matrix %???
differential operators on the circle. 
Using the methods of 
\cite{14}, we show in \cite{15} that a kind of an analogue of the Burnside
theorem holds: $\gcn$ has a unique non-trivial finite irreducible module
(up to going to the contragredient module) and any finite module is a direct 
sum of irreducible modules.

It is a challenging problem to classify all (infinite) simple conformal 
algebras that have a finite module. As we have seen, all of them are 
subalgebras of $\gcn$, some $N$. Incidentally, it follows 
that they have finite 
growth, and it is natural to state a more general problem: classify all
simple conformal algebras of finite growth. 

Note that vertex algebras viewed as conformal algebras  have an exponential
growth.

%6.4. 
\subsection{$q$-analogues \protect\cite{16}} 
\label{sec:6.4}

One can develop a parallel theory of 
$q$-analogues of conformal algebras. They include the ``twisted'' current 
conformal algebras. The general $q$-conformal algebra is that associated 
with the $sin$ algebra \cite{17}.

We hope that this will lead eventually to $q$-analogues of vertex algebras.

\end{document}